\documentstyle[12pt,euscript,leqno]{amsart}
\newtheorem{thm}{Theorem}[section]
\newtheorem{cor}[thm]{Corollary}
\newtheorem{prop}[thm]{Proposition}
\theoremstyle{definition}

\newcommand{\rn}{\Bbb{R}^{\it n}}

\newcommand{\Ho}{H\"{o}lder}
\newcommand{\Gbk}{\overline{\, \Gamma_k}}

\newcommand{\oscsig}{\underset{B_{\sigma}}{\rm osc}}
\newcommand{\oscB}{\underset{B}{\rm osc}}

\newcommand{\sn}{\mathbb{S}^{\it n}}
\newcommand{\sk}{\mathcal{S}_{\it k}}

\newcommand{\pl}{\, + \,}
\newcommand{\mi}{\, - \,}
\newcommand{\Gk}{\Gamma_{\it k}}
\newcommand{\po}{\partial \Omega}

\newcommand{\les}{\, \le \,}

\newcommand{\ges}{\, \ge \,}
\newcommand{\eqs}{\, = \,}
\newcommand{\nd}{\noindent}

\newcommand{\nt}{\notag}

\numberwithin{equation}{section}

\begin{document}

\title[New Maximum Principles for linear elliptic equations]
{New Maximum Principles for linear elliptic equations }
\author{Hung-Ju Kuo${}^\dagger$}
\thanks{${}^\dagger$Research supported by
Taiwan National Science Council}
\address{${}^\dagger$Department of Applied Mathematics, National
Chung-Hsing University, Taichung 402, Taiwan.}
\email{kuohj@@nchu.edu.tw}
\author{Neil S. Trudinger${}^\ddagger$}
\thanks{${}^\ddagger$Research supported by Australian
Research Council Grant.}
\address{${}^\ddagger$Centre for mathematics and Its Applications,
Australian National University, Canberra, ACT 0200, Australia.}
\email{neil.trudinger@@anu.edu.au}

\maketitle
\begin{abstract}
 We prove extensions of the estimates
 of Aleksandrov and Bakel$'$man
 for linear elliptic operators in Euclidean space $\rn$
 to inhomogeneous terms in
 $L^q$ spaces for $q < n$. Our estimates depend on
 restrictions on the ellipticity of the operators determined by certain
 subcones of the positive cone. We also consider some applications
 to local pointwise
 and $L^2$ estimates.

 \end{abstract}

\section{Introduction}

\noindent $\qquad$In this paper, we consider linear second order
partial differential operators $L$ of the form
%
%
\begin{equation} \label{def:L}
 Lu \, := \, a^{ij} D_{ij} u \, ,
\end{equation}
in bounded domains $\Omega$ in Euclidean $n-$space $\rn$. The
operator $L$ is $elliptic$ in $\Omega$ if the coefficient matrix
$\mathcal{A}$ $= [a^{ij}] : \Omega \rightarrow \mathbb{S}^{\it n}$
is positive in $\Omega$. Here $\sn$ denotes the linear space of $n
\times n$ real symmetric matrices and $D^2u = [D_{ij} u] \in \sn$
is the Hessian matrix of second derivatives of an appropriately
smooth function $u: \Omega \rightarrow \Bbb{R}$. The maximum
principle of Aleksandrov and Bakel$'$man
\cite{Alek:60, Alek:68a,  Bak:94}
provides for any solution $u \in C^2(\Omega) \cap
C^0(\overline{\Omega})$ of the inequalities,
%
%
\begin{align} \label{ineq:subsol}
  Lu &\ges -f  \qquad  \; \text{ in } \; \Omega, \\
   u &\les 0    \quad \qquad  \text{ on } \; \partial \Omega, \nt
\end{align}
an estimate
%
%
\begin{equation} \label{prop:maxold}
   \sup_{\Omega} \, u \les C \, \left | \left |
    \frac{f}{\rho_{_n}(\mathcal{A})} \right | \right |_{L^n(\Omega)}
\end{equation}
where $C$ is a constant depending on $n$ and $\Omega$ and the
function $\rho_n$ is given by
%
%
\begin{equation} \label{def:rhoA}
   \rho_{_n} (\mathcal{A}) \eqs  (\det \mathcal{A} )^{\it 1/n} .
\end{equation}
In the special case, where $L$ is the Laplacian, that is
 $\mathcal{A} =  {\it I} $, the exponent $n$ in \eqref{prop:maxold}
 can be improved so that
%
%
\begin{equation} \label{max:newbon}
   \sup_{\Omega} \, u \les C  \parallel f \parallel
   _{L^{^q}(\Omega)}
\end{equation}
for any $q>n/2$, where $C$ is a constant depending on $n, q$ and
$\Omega$. Our concern in this paper is with estimates which lie
between these two extreme cases. As well we shall treat more
precise forms of these estimates, along with
applications to local estimates. \\

 To illustrate the nature of our results we first formulate here an
 extension of the estimates
 \eqref{prop:maxold} \eqref{max:newbon}. The coefficient
 conditions will be expressed in terms of subcones of the positive
 cone in  $\sn$, $\Gamma_n = \{ \mathcal{A} \in \sn \, | \,
 \mathcal{A} > {\rm 0} \}$, determined by the elementary symmetric
 functions $\mathcal{S}_{\it k}$, $ k = 1, \,  \cdots \, n,$ given by
%
%
\begin{equation}  \label{def:Sk}
\mathcal{S}_{\it k} (\lambda) \eqs \sum \lambda _{{\it i}_{\rm 1}} \, \cdots
\, \lambda _{{\it i}_{\it k}}
\end{equation}
for $\lambda = (\lambda_1 , \cdots \, , \, \lambda_n) \in \rn$,
where the summation is taken over all increasing $k$-tuples $\{
i_1 , \cdots \, , \, i_k  \} \subset \{ 1,  \cdots \, , \, n \}$.
Let us first note that for a convex symmetric cone $\Gamma$ in
$\rn$, the dual cone $\Gamma^*$ given by
%
%
\begin{equation}  \label{def:dualcn}
   \Gamma^* \eqs  \left \{ \, \lambda \in \, \rn \, \big | \, \lambda \cdot \mu
   \ges  {\rm 0} \, , \quad \forall \; \mu \in \Gamma \, \right \}
\end{equation}
is closed, convex and symmetric. We associate with the elementary
symmetric function $\sk$  the open cone
%
%
\begin{equation}  \label{def:Gammak}
    \Gamma _k \eqs \left \{ \, \lambda \in \rn \, \big | \;
    \mathcal{S}_{\it j} (\lambda) >  {\rm 0}, \quad {\it j} = {\rm 1,
    \cdots \, , \, k } \, \right \}
\end{equation}

\nd and its closure
%
\begin{equation}  \label{def:Gammab}
 \overline{\Gamma _k} \eqs \left \{ \, \lambda \in \rn \, \big | \;
    \mathcal{S}_{\it j} (\lambda) \ges  {\rm 0}, \quad {\it j} = {\rm 1,
    \cdots \, , \, k } \, \right \}
\end{equation}
which are both convex and symmetric. Clearly, $\Gamma_k \subset
\Gamma_l$ for $k \ge l$ and $\Gamma_1$ is the half-space,
$\Gamma_1 = \{\lambda \in \rn \, | \, \sum \lambda_{\it i} >  {\rm
0} \}$, while $\Gamma_n$ is the positive cone, $\Gamma_n \{\lambda
\in \rn \, | \, \lambda_{\it i} >  {\rm 0} , \, {\it i} = {\rm 1}
, \cdots \, , \,  {\it n} \}$. Note that $\Gamma_k$ can also be
characterized as the component of the positivity set of $\Gamma_k$
which contains $\Gamma_n$, as in \cite{CaNiSp:85}.
Consequently the dual
cones $\Gamma^*_k \subset \Gamma^*_l$ for $k \le l$ with
$\Gamma_1^*$  the closed ray through (1, ... , 1) and $\Gamma_n^*
= \overline{\Gamma_n}$. Corresponding dual functions are
determined as follows. First, we normalize $\sk$ by defining, for
$\lambda \in \Gamma_k$,
%
%
\begin{equation}  \label{def:rhok}
    \rho _ {_k} (\lambda) \eqs
     \left \{ \frac{\sk(\lambda)}{{{\it n} \choose
    {\it k}}
    } \right \}^{1/k}.
\end{equation}
We remark that the function $\rho_{_k}$ is increasing and concave
on the cone $\Gamma_k$, \cite{CaNiSp:85}, and $\rho_{_k} \le
\rho_{_l} $ if $k \ge l$, (Maclaurin inequalities) . The dual
function $\rho^*_{_k}$ is defined on $\Gamma_k^*$ by
%
%
\begin{equation}  \label{def:rhosta}
   \rho^*_{_k} \eqs \inf \left \{ \, {\lambda \cdot \mu \over n } \,
   \big | \, \mu \in \Gamma_k \, , \quad \rho_{_k} (\mu) \ge 1 \,
   \right \}.
\end{equation}
Clearly we have $\rho^*_1(\lambda) = \ell$ where $\lambda=\ell
(1,\cdots , 1)$ and $\rho^*_n(\lambda)=\rho_n(\lambda)= (\prod
\lambda_i)^{1/n}$. As a further example, we may calculate
%
%
%
\begin{align}  \label{form:Gamasto}
   \Gamma^*_2 &\eqs  \left \{\lambda \in \rn \; \big | \;  |\lambda
   | \les {{\rm 1} \over \sqrt{{\it n-}{\rm 1}} } \sum \lambda_{\it i} \right \} , \\
   \rho_2^* &\eqs {1 \over \sqrt{n}} \left \{ \big (\sum \lambda_i \big
   )^2 - (n-1) \, |\lambda|^2 \,
   \right \}^{1/2} . \nt
\end{align}

\nd We shall employ the same notation as above for matrices
 $\mathcal{A} \in \sn$, writing $\mathcal{A} \in \Gk$$(\Gbk , \Gamma^*_k)$
 if the eigenvalues of  $\mathcal{A}, \lambda = \lambda(\mathcal{A}) \in
 \Gk$$(\Gbk , \Gamma^*_k)$ and define $\sk(\mathcal{A})=\sk(\lambda)$,
 $\rho_{_k}(\mathcal{A})$ = $\rho_{_k}(\lambda)$, $\rho_{_k}^*(\mathcal{A})$
 = $\rho_{_k}^*(\lambda)$.
 We can now state the following extension of the estimates \eqref{prop:maxold},
 \eqref{max:newbon}.

%
%
%
%
%
\begin{thm} \label{thm:newmax}
   Let u $\in C^2(\Omega) \cap C^o(\overline{\Omega}) $ satisfy
   \eqref{ineq:subsol} for some coefficient matrix $\mathcal{A}
   \in \Gamma^*_{\it k}$, $1 \le k \le n$, with $\rho^*_k (\mathcal{A}) >
   {\rm 0}$. Then we have the estimate
%
%
\begin{equation} \label{ineq:supnom}
    \sup_{\Omega} \, u  \les C \left | \left | {f \over
   \rho^*_{_k}(\mathcal{A}) } \right | \right | _{L^q(\Omega)}
\end{equation}

\nd for $q = k$ if $k > n/2$ and $q>n/2$ if $k \le n/2$, where $C$
is a constant depending on $n, q$ and $\Omega$.
\end{thm}
\vspace{5mm}

It follows, by approximation as in \cite{GT:83}, that Theorem
\ref{thm:newmax} extends to functions $u \in W^{2,q}_{\ell\,{\rm
oc}}(\Omega) \cap C^0(\overline{\Omega})$. Accordingly we have the
uniqueness result that if $\mathcal{A} \in$ $\Gamma^*_k$,
$\rho^*_k(\mathcal{A})$ $> 0$, $Lu = 0$ a.e. $(\Omega)$, $u = 0$
on $\partial \Omega$, then $u = 0$ in $ \Omega$. Using the example
of Gilbarg and Serrin, (see \cite{GT:83}),

%
%
\begin{equation} \label{examp:Lu}
    Lu \eqs \Delta u \pl \big (\,  -1 \pl {n-1 \over 1 - \alpha } \,
    \big ) \, { x_i x_j \over |x|^2}\; D_{ij} u\, ,
    \quad  \alpha < 1
\end{equation}

\nd with solution $u$ given by

%
%
\begin{equation} \label{funct:uabs}
     u(x) \eqs \begin{cases}|x|^{\alpha}&, \quad  \text{ if } \alpha \ne 0 \\
       \log |x|&, \quad \text{ if } \alpha = 0.
       \end{cases}
\end{equation}

\nd satisfying $ u \in W^{2,q}_{\ell\,{\rm oc}}(\rn)$ if and only
if $q <{n \over 2-\alpha}$, we infer that the exponent $q$ in
Theorem \ref{thm:newmax} cannot be improved (i.e. $q$ can not be smaller than $k$).
To see this, we note,
for example from \cite{TW:99}, (see also \cite{Tr:06}), that
$\mu =(\mu_1\, , \, \cdots \, , \, \mu_n) \in
\Gamma _k$ implies

%
%
     \begin{equation} \label{ineq:mui}
         k(n-1) \, \mu_i \pl (n-k) \, \sum_{j\ne i} \mu_j \ges 0
     \end{equation}
\nd for any $i=1, \cdots \, , \, n$. For the coefficient matrix
$\mathcal{A}$ in \eqref{examp:Lu}, we then have
   \begin{equation*}
       \lambda(\mathcal{A} )
       \eqs \big ({\rm 1}, \, \cdots \, , \,
       {\rm 1}, \, {{\it n}-{\rm 1} \over {\rm 1}
        - \alpha} \big ) \, \in \, \Gamma_{\it k}^* \, ,
   \end{equation*}
\nd provided $\alpha \le 2- n/k$. Consequently for
${n \over 2} < k < n, \, q < k$, we choose $\alpha = 2 - n/k$
to get a counterexample. The case $k=n$,
follows by a slight modification,
taking $2-n/q < \alpha < 1$ for $q < n$. For the case $k \le n/2$, if $ q<n/2$,
we get a counterexample
with $\alpha < 0$ while for $q =n/2$,
we may modify
\eqref{funct:uabs} by taking $\alpha = 0$ and for $\varepsilon > 0$,

  %
  %
  \begin{equation} \label{def:upsln}
     u_{\varepsilon} (x)\eqs \begin{cases} \log |x|&, \quad  \text{ for } |x| \ges \varepsilon \\
                {1 \over 2}  \left ( {|x|^2 \over \varepsilon^2 } -1 \right )
     + \log\varepsilon &, \quad \text{ for } |x| < \varepsilon.
                                 \end{cases}
  \end{equation}
 \nd In this case $\|Lu_{\epsilon} \|_{L^q(\Omega)}$ is uniformly bounded in
  $\varepsilon$ but $\inf u_{\varepsilon}
   \rightarrow - \infty$ as $\varepsilon \rightarrow 0$.
In this connection we mention the recent work of
Astala, Iwaniec and Martin, \cite{AsTwMa:06}, for $n=2$, where an
estimate of the form
%
%
\begin{equation} \label{bound:supuf}
   \sup_{\Omega} \, u \les C \, \| f\|_{L^q(\Omega)} \, ,
\end{equation}

\nd is derived for solutions of \eqref{ineq:subsol} provided $q >
{2K \over K+1}$, where  $K = \sup_{\Omega} \lambda_{\rm
max}(\mathcal{A})/ \lambda_{\rm min}(\mathcal{A})$ denotes the
ellipticity constant of $\mathcal{A}$. The operator
\eqref{examp:Lu} may also be used to show that their estimate
\eqref{bound:supuf} is sharp, \cite{AsTwMa:06}. \\

 In the next section we will in fact prove a stronger version
of Theorem \ref{thm:newmax}, where the $L^q$ norm is taken over
the upper $k-$contact set of $u$ in $\Omega$. In the following
section, we will consider sharp versions of the estimate
\eqref{ineq:supnom} in the cases $k>n/2$, using the Greens
function from \cite{TW:97}. Finally in Section 4, we prove a
corresponding local maximum principle and indicate the relevant
extensions of other local estimates such as the Harnack and Holder
estimates, \cite{GT:83, KrSa:81}. As an application of the local
maximum principle, we obtain an extension of \eqref{ineq:supnom}
for uniformly elliptic operators, with the constant $C$
depending only on $n,k$ and $|\Omega|$, analogous to \cite{Cab:96}.

Some of this paper, in particular Theorem \ref{thm:newmax},
was proved by the second author several years ago
and presented at
various meetings. The two authors have also obtained discrete
versions of the case
$k=n$, ( \cite{KT:96, KT:pucci, KT:00}).
It would also be interesting to have
corresponding discrete versions of the estimates in this paper.

\vspace{10mm}

%
%
\section{Reduction to Hessian Equations}

\vspace{4mm}

 For $k=1, \cdots , n$, the $k$-Hessian operator $F_k$ is defined
 on $C^2(\Omega)$ by
%
%
\begin{align} \label{def:FkSk}
 F_k[u] &\eqs
    \sk({\it D}^{\rm 2}{\it u}) \\  &\eqs \big [ {\it D}^{\rm 2}{\it u} \big ]_{\it
    k} \nt
\end{align}

\nd where for an $n \times n$ real matrix, $\mathcal{A}\;$,
$[\mathcal{A}]_{\it k}\;$ denotes the sum of its $k \times k$
principal minors. The operator $F_k$ is related to the linear
operator in \eqref{def:L} through the following proposition.
%
%
\begin{prop} \label{prop:innerp}
  For any matrices $\mathcal{A} \in \Gk$, $\mathcal{B} \in
  \Gamma^*_{\it k}$, $k=1, \cdots \, , \, n$, we have the
  inequality,
  %
  %
\begin{equation} \label{ineq:ABiner}
   \rho_k(\mathcal{A}) \rho^*_{\it k}(\mathcal{B}) \les {{\rm 1}
   \over {\it n} } \, \mathcal{A} \cdot \mathcal{B} \, .
\end{equation}
\end{prop}
\nd {\bf Proof. } If we fix the matrix $\mathcal{B}$ = $[b_{ij}]$
and minimize the inner product $\mathcal{A}\cdot \mathcal{B}$
 on the set where $\mathcal{A}$$\cdot
\mathcal{B} \ge$ 0, $\sk(\mathcal{A})$ = 1, we obtain at a
critical point $\mathcal{A}$,
$$
  \mathcal{B} \eqs {\it c} \,  D \sk (\mathcal{A})
$$
for some constant $c$. Hence with respect to an orthonormal basis
of eigenvectors of $\mathcal{A}$, the matrix $\mathcal{B}$ is also
diagonal with eigenvalues $\mu = (\mu_1, \, \cdots \, , \, \mu_n)$
given by
$$
  \mu_i \eqs c \, D_i\sk(\lambda)\, ,  \qquad  {\it i} = {\rm 1,}
  \cdots \, , \, {\it n} ,
$$
\nd where $\lambda = \lambda(\mathcal{A})$, \cite{CaNiSp:85}. The
inequality \eqref{ineq:ABiner} then follows from the definition
\eqref{def:rhosta}. \hspace{100mm}
\newline $\Box$

Let $u \in C^2(\Omega)$ satisfy the differential inequality
\eqref{ineq:subsol}. From Proposition \ref{prop:innerp}, we then
have
%
%
\begin{align} \label{ineqs:rhoDDu}
    \rho_{_k}(-D^2u) \rho^*_{_k} (\mathcal{A}) &\les  - \, {{\rm 1} \over {\it n}} \, {\it Lu} \\
           &\les {1 \over n} \,  f \nt
\end{align}

\nd where $-D^2u \in \overline{\Gamma_k}\,$, $\mathcal{A}$ $\in
\Gamma^*_k$ . Consequently, replacing $u$ by $-u$ we have the
differential inequality,
%
%
\begin{equation} \label{ineq:Fkpsi}
   F_k[u] \les \psi \, , \qquad \qquad \text{ where } \qquad
   \psi  \eqs {n \choose k} \,
   \left ( {f \over n \, \rho^*_{_k} (\mathcal{A}) } \right )^k
\end{equation}

\nd holding on the subset of $\Omega$ where $D^2u \in
\overline{\Gamma_k}$ , that is where the function $u$ is
$k$-convex. The estimate \eqref{ineq:supnom} is accordingly
reduced to the existence and estimation of solutions of Hessian
equations, with inhomogeneous terms in $L^{^p}$ spaces. Indeed if
$u \in C^2(\overline{\Omega})$ is a $k-$convex function on
$\Omega$, satisfying \eqref{ineq:Fkpsi}, with $u=0$ on $\partial
\Omega$, then it follows readily from the Wang Sobolev
inequality, \cite{ChW:01, Wa:94}, using Moser iteration,
 that

%
%
\begin{equation} \label{ineq:supulp}
   \sup_{\Omega} \, u\les C \, \| \psi \|^{1/k} _{_{L^{^p}(\Omega)}}
\end{equation}

\noindent where $p=1$ for $k > {n \over 2}$, $ p > {n \over 2k}$
for $k \les {n \over 2}$ and $C$ depends on $k, \, n, \,  p, $ and
diam$\Omega$. \\

From the estimate \eqref{ineq:supulp}, we can prove Theorem
\ref{thm:newmax}, as stated, through an existence theorem for
Hessian equations. Let $\Omega_o$ be a uniformly ($k-1$)-convex
domain, containing $\Omega$, with boundary $\po_o \in C^{\infty}$
and set

%
%
\begin{equation} \label{def:psiprm}
    \psi' \eqs F_k[u] \chi_{_{\Omega_k}} \; ,
\end{equation}

\nd where $\Omega_k  = \Omega^-_k$ is the {\it lower $k-$contact set} of
$u$ in $\Omega$ given by

%
%
\begin{equation} \label{def:lconst}
   \Omega^-_k \eqs \big \{ \, x \in \Omega \, \big | \,
   \exists \; k-\text{convex} \; v \in C^2(\Omega) \text{ satisfying }
   v \le u \text{ in } \Omega , \; v(x) = u(x) \; \big \}.
\end{equation}

\nd Clearly for any $x \in \Omega^-_k , \, D^2u(x) \in \overline{\Gamma_k}$.
By replacing $\Omega$ if necessary by
a strictly contained subdomain we may assume $u \in
C^2(\overline{\, \Omega})$, and $\psi '  \in C^o( \overline{\,
\Omega})$. For $\psi_o \in C^{\infty} (\overline{\Omega})$
satisfying $\psi_o
>0$ in $\overline{\, \Omega}_o$, $\psi' < \psi_o$ in $\Omega$,
we define $u_o \in C^{\infty}(\overline{\,\Omega}_o)$ to be the
unique $k$-convex solution of the Dirichlet problem,
%
%
%
\begin{align} \label{Problm:Fk}
     F_k[u_o] &\eqs \psi_o \qquad \qquad  \text{in } \; \Omega_o \quad , \\
     u_o &\eqs 0 \qquad \qquad \; \text{on } \;  \po_o \; \, . \nt
\end{align}
The existence of $u_o$ is guaranteed by the existence theorem of
Caffarelli, Nirenberg and Spruck [8]; (\,see also \cite{Tr:95b}\,).
We claim that a comparison principle holds namely $u \ge u_0$ in
$\Omega$. To see this we suppose there exists $y \in \Omega$ such that

$$
  u_0(y) \, - \, u(y) \eqs \, \sup_{\Omega} (u_0 \, - \, u) \, > \, 0.
$$

\nd Since $u_0$ is $k-$convex, we must have $y \in \Omega^-_k$.
But then $D^2u_0(y) \les D^2u(y)$ implies $F_k[u_0] \les F_k[u]$,
which contradicts \eqref{Problm:Fk}.
By letting $\psi_0$
approach to $\psi '$, we then obtain the estimate \eqref{ineq:supulp}, with
$\psi = \psi ', \Omega = \Omega _0$.
Hence letting $\Omega^+_k (u) = \Omega^-_k (-u) $
  denote the {\it upper k-contact set} of
  $u$ in $\Omega$, we obtain the estimate

%
%
\begin{equation} \label{est:supu}
   \sup_{\Omega} u \les C \left | \left | {f \over
   \rho^*_{_k}(\mathcal{A}) } \right | \right | _{L^q(\Omega^+_k)},
\end{equation}

\nd which is a more precise version of \eqref{ineq:supnom}. In the
next section, we shall provide a proof
of \eqref{ineq:supulp} in the cases $k > n/2$,
using \cite{TW:97}, which leads to sharper
versions of \eqref{est:supu}.

\bigskip

\bigskip

\bigskip

\bigskip
%
%
\section{Refinements}

 In this section we refine the estimate \eqref{est:supu}
 by using an argument analogous to that of Aleksandrov
 and Bakel$'$man  for the case $k=n$; see \cite{GT:83}. The fundamental idea is to
 replace cones by the graphs of the Green's functions in the cases $k>n/2$.
 For this and the following section,
 we need some aspects of the theory of Hessian
 measures developed by Trudinger and Wang in
\cite{TW:97,  TW:99, TW:02}. First we recall the general
definition of $k-$convexity. Namely an upper semi-continuous
function $u : \Omega \rightarrow [- \infty, \, \infty )$ is called
$k-convex$ in $\Omega$ if any quadratic polynomial $p$ for which
$u-p$ has a local maximum in $\Omega$, satisfies $F_k[p] \ges 0$.
General $k-$convex functions may be approximated by smooth ones
through mollification. Indeed let us define $\Phi^k(\Omega)$ to be
the set of proper $k-$convex function, that is those
$\equiv$\hspace{-2.6mm}/ $- \infty$
 on any component of $\Omega$. Then $\Phi^k(\Omega) \subset \Phi_1(\Omega)
 \subset L^1_{\ell\,{\rm
oc}}(\Omega)$ and the mollification $u_h$ of
 $u \in  \Phi^k(\Omega)$ satisfies $u_h \downarrow u$ as $h \rightarrow 0$,
 $u_h \in \Phi^k (\Omega')$ for any
 $h < \text{dist}(\Omega', \partial \Omega)$, $\Omega' \subset \subset \Omega$.
 The main result of \cite{TW:97,  TW:99} is that for any $u \in \Phi^k(\Omega)$, there
 exists a Borel measure $\mu_k[u]$ such that

 %
 %
\begin{align} \label{def:muk1}
 &\quad (i)  \quad \mu_k[u](e) \eqs \int_{e} F_k[u] \qquad
      \text{ for any } u \in \Phi^k(\Omega) \cap C^2(\Omega),
       \qquad \qquad \qquad \qquad\\
&\text{and} \nt \\
&\quad (ii)  \quad \text{if } u_m \rightarrow u \text{ a.e. in } \Omega,
 u_m, u \in \Phi^k(\Omega), \nt \\
 &\qquad \qquad \qquad  \text{ then }
  \mu_k[u_m] \rightarrow \mu_k[u] \text{ weakly as measures }. \nt
 \end{align}

\nd The case
 $k > n/2$ is proved in \cite{TW:97}. Here $\Phi^k(\Omega) \subset
 C^{0,2-n/k} (\Omega)$ and a.e. convergence is equivalent to uniform convergence.
 From the weak continuity ($ii$) of \eqref{def:muk1}, it follows that
 there exists a unique Green's function for $F_k$. That is for any
 point $y \in \Omega$, there exists a function $G_y \in \Phi^k(\Omega)$
 such that
  %
 %
\begin{align} \label{def:greenk}
  \mu_k[G_y] &\eqs \delta_y \; ,  \qquad \qquad \\
  G_y &\, \longrightarrow \;\; 0 \; \text{ on } \partial \Omega. \nt
\end{align}

\nd Here $\delta_y$ denotes the Dirac delta measure at $y$ and we
also need to assume that $\partial \Omega \in C^2$ is uniformly
$(k-1)-convex$, that is the principal  curvatures $(\kappa_1 \, ,
\, \cdots \, , \,  \kappa_{n-1})$ of $\Omega$ lie in the cone
$\Gamma_{k-1}$ in $\mathbb{R}$$^{n-1}$. When $k>n/2, G_y \in
C^{0,2-n/k}(\overline{\Omega}\,)$. The uniqueness is more
difficult in the cases $k\les n/2$, \cite{TW:02}.

Moreover, from the interior gradient bound  \cite{Tr:97}, we
always have $G_y \in C^{0,1}(\overline{\Omega} - \{y\})$. It is
easy to show, for example by smoothing the cusp, that the Green's
function $G_y$ for a ball $B_R(y)$ of center $y$ and radius $R$ is
given by

 %
 %
\begin{equation} \label{def:Gy}
G_y(x) \eqs \left \{
\begin{aligned}
&{1 \over (2-{n \over k}) \left [  \binom{n}{k}  \omega_n \right
]^{1/k}} \left ( |x-y| ^{2-n/k} - R^{2-n/k} \right ) \quad, \text{
if }
    k \ne {n \over  2}  ,  \\
&{1 \over \left [ \binom{n}{k} \omega_n \right ]^{1/k}} \; \log
|x-y| \qquad \qquad \qquad \qquad \qquad \;
 ,\text{ if } k = {n \over  2}.
\end{aligned}
\right.
\end{equation}

\nd In the case $k=n$ and convex $\Omega$, the Green's function
 $G_y$ is the function whose
graph is a cone with vertex at $(y, G_y(y))$ and base $\partial
\Omega$. For general $\Omega$ and  $k>n/2$, $G_y$ will have a cusp-like behavior at $y$,
as exemplified by \eqref{def:Gy}. It is also shown in
\cite{TW:97}, that the monotonicity property of the
Monge-Amp\`{e}re measure in  the case $k=n$, extends to Hessian
measures for $k \les n$. Namely if $u, v \in \Phi^k(\Omega)$
satisfy $u \les v$ in $\Omega$, $u=v$ continuously on $\partial
\Omega$, then $\mu_k[u](\Omega) \ge \mu_k[v](\Omega)$. From this we
also have a comparison principle, namely if $u, v \in
\Phi^k(\Omega), u \le v$ continuously on $\partial \Omega$,
$\mu_k[u] \ge \mu_k[v]$ in $\Omega$, then $u \le v$ in $\Omega$. As
a consequence we obtain an estimate for the Greens' function, in
the cases $k>n/2$, by comparison with \eqref{def:Gy}. That is

%
%
\begin{equation} \label{inf:Gy}
  \inf_{\Omega} G_y \eqs G_y(y) \ges \;
 \; - \,\; { (\text{diam } \Omega )^{2-n/k} \over
   (2- {n \over k}) \left [ \binom{n}{k} \omega_n \right ]^{1/k} }
\end{equation}

\nd Now returning to the proof of Theorem \ref{thm:newmax}
in Section 2, we let $u \in \Phi^k(\Omega)\cap
C^2(\Omega) \cap C^0(\overline{\Omega})$ satisfy

%
%
\begin{align} \label{ineq:Fk}
 F_k[u] &\les \psi  \qquad \text{ in } \;  \Omega \\
 u &\eqs 0   \qquad \text{ on } \;  \partial \Omega \nt
\end{align}

\nd where $\partial \Omega$ is $(k-1)-$convex and $k>n/2$. Then
for any point $y \in \Omega$, we have
%
%
\begin{equation} \label{ineq:uGy}
  u \les \left [ {u(y) \over G_y(y) } \right ] \, G_y
\end{equation}

\nd since $\mu_k [G_y] = 0$ in $\Omega - \{y\} $.
Hence, by the monotonicity property of $\mu_k$ and its $k-$homogeneity,
we obtain

\begin{align*}
  \left [ {u(y) \over G_y(y) }\right ]^k
  &\eqs  \left [ {u(y) \over G_y(y) }\right ]^k \mu_k[G_y](\Omega) \\
  & \les \mu_k[u] (\Omega) \\
  &\les \int _{\Omega} \psi
\end{align*}

\nd so that we have the precise estimate

%
%
\begin{align} \label{ineq:uGy2}
 -u(y) &\les -G_y(y) \left ( \int_{\Omega} \psi \right )^{1/k} \\
 &\les  {(\text{diam} \, \Omega)^{2-n/k} \over (2-n/k)
 \left [\binom{n}{k} \omega_n \right ]^{1/k}}
 \left ( \int_{\Omega} \psi \right )^{1/k} \nt
\end{align}

\nd by virtue of \eqref{inf:Gy}. Accordingly we obtain
\eqref{est:supu} with constant $C$ given by

%
%
\begin{equation} \label{const:C}
  C \eqs {(\text{diam} \, \Omega)^{2-n/k}
  \over n(2-n/k)(\omega_n)^{1/k}}.
\end{equation}

\nd Instead of using the Green's function  $G_y$, we may
use the function
%
%
\begin{equation} \label{def,Wy}
   w_y \eqs - \, {G_y \over G_y(y) }
\end{equation}

\nd which can be obtained independently as the weak solution of the homogeneous
Dirichlet problem,

%
%
\begin{align} \label{funct,Fk}
   F_k[w_y] &\eqs 0 \qquad \text{ in }  \; \Omega - \{y\} , \\
    w_y &\eqs 0 \qquad \text{ on } \partial \Omega \, , \nt \\
    w_y(y)&\eqs 1 \, , \nt
\end{align}

\nd for example by using the Perron process. It then follows directly
from \eqref{funct,Fk} that $w_y \in \Phi^k(\Omega)
\cap C^{2-{n \over k}} (\overline{\Omega})
\cap C^{0,1}(\overline{\Omega}-\{y\})$ and
moreover in the estimate \eqref{ineq:uGy2},

%
%
\begin{equation} \label{ineq:uy}
   -\, u(y) \les \left \{ {1 \over \mu_k[w_y] }
   \int_{\Omega} \psi \right \} ^{1 \over k}.
\end{equation}

The quantity $\mu_k[w_y]$ is an extension to $n/2 < k \le n$
of the volume of the polar of $\Omega$, with respect to $y$,
in the case $k=n$. The best constant $C$ in \eqref{est:supu}
is thus given by

%
%
\begin{equation} \label{Const:C}
   C \eqs {1 \over n}  \binom{n}{k}  ^{1 \over k}
    \, \sup _{y  \in \Omega}
   \Big \{ \mu_k [w_y] \Big\}^{-\, {1 \over k}}.
\end{equation}

\nd
Note that the cruder estimate \eqref{const:C}
may be proved directly from
\eqref{def:Gy} using the comparison principle.

\bigskip

\bigskip

\bigskip

\bigskip
%
%
\section{Local estimates}

\bigskip

In this section, we consider the full linear operator,
%
%
\begin{equation} \label{def:Lulst}
  Lu \, := \, a^{ij} D_{ij} u \pl b^iD_iu \pl cu
\end{equation}

\nd under the hypothesis, $\mathcal{A}$ = $[a^{ij}] \in \Gamma^*_k$,
with

%
%
\begin{align} \label{condit:rhostr}
   \rho^*_k(\mathcal{A}) &\ges  \rho_{_0} \\
   |\mathcal{A}| &\les a_0 \nt
\end{align}

\nd where $\rho_{_0}$ and $a_0$ are  positive constants and $b^i,
\, c  \in L^{\infty}(\Omega)$. Note that by \eqref{def:rhosta},
condition \eqref{condit:rhostr} implies $L$ is uniformly elliptic, as on
$\Gamma^*_k$,

%
%
\begin{align} \label{ineq:lambda}
   \lambda_n &\ges  \lambda_1^{1-n} \rho^n_n(\mathcal{A}) \, , \\
      &\ges a_0^{1-n} \big ( \rho^*_k(\mathcal{A}) \big )^{\it n} \, , \nt \\
        &\ges   a_0^{1-n} \rho^n_0  \, , \nt
\end{align}

\vspace{2mm}

\nd where $\lambda_1 \, ,  \, \lambda_2 \, ,  \, \cdots \, , \, \lambda_n$
denote the eigenvalues of $\mathcal{A}$ in decreasing order. Local pointwise
estimates for $L$ then follow from Theorem \ref{thm:newmax} as with the
uniformly elliptic case $k=n$. We first consider an extension
of the local maximum principle of Trudinger \cite{GT:83, Tr:80}.
%
%
\begin{thm} \label{thm:locmax}
 Let u $\in C^2(\Omega)\cap C^0(\overline{\Omega})$ satisfy
%
%
\begin{align} \label{equatin:LuB}
   Lu &\ges - f \qquad \qquad {\rm in } \; \; \Omega \cap B \, , \\
   u &\les  0    \quad \qquad \qquad  {\rm on } \;   \;
     \partial (\Omega \cap B ) \, , \nt
\end{align}

\nd for some ball $B = B_R(y) \subset \rn$ , $f \in
L^q(\Omega) $ where $q = k$ for $k > n/2$ , $q > n/2$ if
$k \les n/2$. Then for any concentric ball
$B_{\sigma} = B_{\sigma R}(y), \; 0 < \sigma < 1$, and $ \, p> 0$,
we have the estimate
%
%
\begin{equation} \label{Estimt:Locmax}
  \sup_{\Omega \cap B_{\sigma}} \, u \, \les C \,
   \left \{ \left ( R^{-n} \int _{\Omega \cap B} (u^+)^p \right )
   ^{1/p}  \pl {R^{2-n/q}
    \over \rho_{_0} } \| f \| _{L^q(\Omega \cap B)} \right \}
\end{equation}

\nd where C is a constant depending on $\sigma, p , n $ and
$a_0/\rho_{_0}$, $\sup |b|R/\rho_{_0} \, , \,  \sup |c|
R^2/\rho_{_0}$
\end{thm}
\begin{pf}
 For the proof of Theorem \ref{thm:locmax},
 we cannot directly employ the proof  in \cite{GT:83,Tr:80}
 but instead we may use that given in \cite{Tr:95a},
 which we indicate briefly here. First, in view of the scaling
 $x \rightarrow x/R, f \rightarrow f/R^2$, we may take $R=1, y=0$.
 Let

%
%
\begin{equation} \label{def:etacut}
   \eta \eqs \left [ (1 - |x|^2)^+ \right ] ^{\beta}
\end{equation}

\nd for $\beta \ge 1$, to be chosen. Setting $v = \eta (u^+)^2$,
we compute in $\Omega \cap B \cap \{u>0\}$

%
%
\begin{align} \label{bound:Lov}
    L_0v &\, := \, a^{ij}D_{ij} v \\
    &\eqs (u^+) ^2  a^{ij}D_{ij} \eta \pl 2a^{ij}D_i \eta
    D_j u^2  \pl \eta a^{ij}D_{ij} u^2 \nt \\
    &\ges  - \, C \, ( |\mathcal{A}| {\it + |b| + |c|})
     \eta ^{\it 1-2/\beta} {\it u^2} \pl {\it 2 \eta u f}  \nt
\end{align}

\nd where $C$ depends on $n, \lambda_1/\lambda_n$ and $\beta$. Now
following \cite{Tr:95a}, we apply Theorem \ref{thm:newmax}, to obtain
%
%
\begin{equation} \label{ineq:Locbal}
  \sup_{\Omega \cap B_{\sigma}} \, v \, \les C \,
   \left \{
   \left \| v^{1-2/\beta } (u^+)^{4/\beta}
     \right \|
   _{L^q(B\cup \Omega)}  \pl {1 \over \rho_{_0}} \;
   \| v^{1/2} f \| _{L^q(B\cup \Omega)}
    \right \} ,
\end{equation}

\nd from which we deduce \eqref{Estimt:Locmax}, by
taking $\beta = 4q/p$.
\end{pf}

We remark that we may only assume $b \in L^{2q}(\Omega)$,
$c \in L^q(\Omega)$ in Theorem \ref{thm:locmax}. Also, as
remarked in \cite{Tr:95a} for the case $k=n$, we obtain by
choosing $R$ sufficiently large in Theorem \ref{thm:locmax},
the following variant of the Aleksandrov-Bakel$'$man
principle for uniformly elliptic operators.
%
%
\begin{cor} \label{coll:max}
 Let u $\in C^2(\Omega)\cap C^0(\overline{\Omega})$ satisfy
\eqref{def:L},
\eqref{ineq:subsol} under the above hypothesis on $L$.
Then we have the estimate

%
%
\begin{equation} \label{bound:supu}
   \sup _{\Omega} u \les C \, |\Omega | ^{2/n-1/q} \left \|
    {f \over \rho^*_k(\mathcal{A}) } \right \|_{L^q(\Omega)},
\end{equation}

\nd where $q=k$ if $k > n/2, \, q > n/2 $ if $k \les n/2$ and
$C $ is a constant depending on $n, q, \,
a_{_0}/\rho_{_0}$.

\end{cor}

The special case $k=n$ of \eqref{bound:supu} was found
differently by Cabr\'{e} in \cite{Cab:96}.

  For solutions, the {\it \Ho} and {\it Harnack}  estimates
of Krylov and Safonov, (see \cite {GT:83, KrSa:81, Tr:80, Tr:95a},
extend automatically to inhomogeneous terms in lower $L^p$ spaces.
This is readily seen, for example by following the proof in
\cite {GT:83, Tr:80}.

 %
 %
 \begin{thm} \label{thm:Harnak}
  Let u $\in C^2(\Omega)$ satisfy $Lu = f$ in $B = B_R(y) \subset \Omega$.
  Then for any concentric ball $B_{\sigma} =B_{\sigma R}(y)$,
  $0 < \sigma < 1$,  we have the estimate

 %
 %
 \begin{equation} \label{ineq:maxmin}
  \oscsig   u  \les C \, \sigma ^{\alpha} \left \{ \oscB u \pl R^{2-n/q}
  \left \| \, {f \over \rho^*_k(\mathcal{A}) }  \, \right \|_{L^q(B)} \right \}
\end{equation}
 \end{thm}

\nd where $q$ is as in Theorem \ref{coll:max},
$\alpha > 0$ depends on $n, a_0/\rho_0 $ and $C$ depends
on $n, q, a_0 / \rho_0$, $\sup |b|R / \rho_0 \, , \,
 \sup |c|R^2  \rho_0$. Furthermore if $u \ge 0$ in $B$, then
 for any $0 < \sigma \, , \, \tau < 1$,
 %
 %
 \begin{equation} \label{ineq:supinf}
   \sup_{B_\sigma} \, u \les C \,
 \left \{ \inf_{B_{\tau}} u \pl R^{2-n/q}
  \left \| \, {f \over \rho^*_k(\mathcal{A} )}  \, \right \|_{L^q(B)} \right \}
\end{equation}

\nd where $C$ depends on the same quantities as in
\eqref{ineq:maxmin} together with $\sigma$ and $\tau$.

We remark that as in Theorem \ref{thm:locmax}  we need only assume
$b / \rho_0 \in L^{2q} (B), c / \rho_0 \in L^q(B)$.
Other local estimates which depend on the Aleksandrov-Bakel$'$man maximum
principle  also extend in a corresponding way.

Finally we remark on an interesting relationship between the case $k=2$ and second
derivative estimates. Indeed we first note another characterization of $\Gamma^*_2$,
namely
%
%
\begin{equation} \label{formua:Gamma}
  \Gamma^*_2 \eqs \left \{ \, \lambda  \in  \rn \, \Big | \,  \left \|
{({\it n}-{\rm 1}) \over {\rm tr} \mathcal{A}} \, \mathcal{A} \mi
 {\it I}
\right \|_{\rm 2} \,  < \, {\rm 1}
\right \}
\end{equation}

\nd Consequently, we have by perturbation from the case $L = \Delta$,
(see \cite{GT:83}),
that if $u \in C^2(\Omega)$, $Lu = f $ in $ \Omega$  and  \eqref{condit:rhostr}
holds for $k = 2$ then for any
$\Omega' \subset \subset \Omega$,

%
%
\begin{equation} \label{ineq:secord}
\left \| D^2 u \right \|_{L^2(\Omega ')} \les C \, \left \{
\| u \|_{L^2(\Omega)}
 \pl
\left \| {f \over \rho^*_2(\mathcal{A}) } \right \| _{{\it L}^{\rm 2}
(\Omega)}
 \right \}
\end{equation}

\nd where $C$ depends on $n, \, \Omega, \, \Omega ' , \, a_0 / \rho_{_0} , \,
b_0 / \rho _{_0}$ and $c_0 / \rho_{_0}$. Now if we take
$n = 3$ and apply Corollary \ref{condit:rhostr}, we obtain the full estimate

%
%
\begin{equation} \label{bound:wtoto}
  \| u \| _{_{W^{2,2}(\Omega')}} \les C \, \left \|
  {f \over \rho^*_{_2}(\mathcal{A})}
  \right \|
  _{{\it L}^{\rm 2} (\Omega)}.
\end{equation}

\nd if $ L$ is of the form \eqref{def:L}
and $ u \in C^2(\Omega)\cap C^0(\overline{\Omega})$
vanishes on $\partial \Omega$. If $\partial \Omega  \in C^{1,1}$,
then we may replace $\Omega '$ by $\Omega$.

Note that the estimates of this section
also extend by approximation
to functions $u$ in  Sobolev spaces $W_{\rm loc}^{2,q}(\Omega)$.

\vspace{20mm}


\bigskip

\bigskip

\bigskip

\bigskip

\end{document}